\documentclass[12pt]{article}
\usepackage[dvips]{epsfig}
\usepackage{amssymb}
\parskip=\medskipamount

\begin{document}

\def\pmb#1{\setbox0=\hbox{#1}%
  \kern-.015em\copy0\kern-\wd0
  \kern.05em\copy0\kern-\wd0
  \kern-.015em\raise.0433em\box0 }

\def\bdJ{\pmb{$J$}}  
\def\bdw{\pmb{$w$}}  
\def\bds{\pmb{$s$}}  
\def\bdF{\pmb{$F$}}  
\def\bd0{\pmb{$0$}}  
\def\bdf{\pmb{$f$}}
\def\bdy{\pmb{$y$}}
\def\bdx{\pmb{$x$}}
\def\bdz{\pmb{$z$}}
\def\bdu{\pmb{$u$}}
\def\bdA{\pmb{$A$}}
\def\bda{\pmb{$a$}}
\def\caL{\cal L}
\def\N{\mathbb N}
\def\R{\mathbb R}
\def\C{\mathbb C}
\def\bdab{\bar{\pmb{$a$}}}
\def\bdom{\pmb{$\omega$}}
\def\be{\begin{equation}}
\def\ee{\end{equation}}
\def\bea{\begin{eqnarray}}
\def\eea{\end{eqnarray}}
\def\beas{\begin{eqnarray*}}
\def\arccos{{arccos}}
\def\eeas{\end{eqnarray*}}
\def\ds{\displaystyle}
\def\b{\phantom{-}}
\newtheorem{dfn}{Definition}
\newtheorem{thm}{Theorem}
\newtheorem{lem}{Lemma}
\newtheorem{cor}[thm]{Corollary}

\renewcommand{\labelenumi}{\Roman{enumi}.}

\date{}
\addtolength{\baselineskip}{0.3\baselineskip}
\renewcommand{\theequation}{\arabic{section}.\arabic{equation}}
\title{\Large\bf Phase properties of exponentially-fitted symmetric multistep methods for $y''=f(x,y)$}
\author{Hans Van de Vyver\\
{Department of Mathematics, Katholieke Universiteit Leuven,}\\
{Celestijnenlaan 200 B, B-3001 Heverlee, Belgium}\\
e-mail: hans$\_$vandevyver@hotmail.com\\
Tel.: +32-16-327048;\qquad Fax: +32-16-327998}
\maketitle
\begin{abstract}
A convenient tool to obtain numerical methods specially tuned on oscillating functions is exponential fitting. Such methods are needed in various branches of natural sciences, particularly in physics, since a lot of physical phenomena exhibit a pronounced oscillatory behavior. Many exponentially-fitted (EF) symmetric multistep methods for $y''=f(x,y)$ are already developed. To have an idea of the accuracy we examine their phase properties. The remarkably simple expression of the phase-lag error obtained in Theorem~\ref{thmpl} allows to draw quantitative conclusions on the merits of each EF version. 
\end{abstract}

\vspace*{1.5cm}
\noindent PACS Subject Classification~: 02.60; 02.70.Bf; 95.10.Ce; 95.10.Eg; 95.75.Pq.

\vspace*{0.5cm}

\vspace*{0.5cm}

\newpage
\setcounter{equation}{0}
\section{Introduction}
Significant efforts were undertaken over the years to promote {\em linear multistep methods} 
\be
\ds\sum_{j=0}^J \alpha_j\,y_{n+j}=h^2\,\ds\sum_{j=0}^J \beta_j\,f_{n+j},\qquad n=0,1,\ldots,
\label{multi}
\ee
where $y_{n+j}$ is an approximation to $y(x_{n+j})$ and $f_{n+j}=f(x_{n+j},y_{n+j})$, as highly competitive solvers for the special second-order initial value problem
\be
y''=f(x,y),\qquad y(x_0)=y_0,\qquad y'(x_0)=y_0'.
\label{ODE2}
\ee
The method (\ref{multi}) is characterized by the polynomials $\rho$ and $\sigma$ where
\[
\rho(\zeta)=\ds\sum_{j=0}^J\alpha_j\,\zeta^j,\qquad \sigma(\zeta)=\ds\sum_{j=0}^J\beta_j\,\zeta^j,\qquad \zeta\in\C.
\]
We associate with method (\ref{multi}) the following linear functional
\be
{\mathcal L}[h,{\bf a}]z(x)=\ds\sum_{j=0}^J \alpha_j\,z(x+j\,h)-h^2\,\ds\sum_{j=0}^J \beta_j\,z''(x+j\,h),
\label{func}
\ee
where ${\bf a}$ is the vector of the coefficients ${\bf a}=[\alpha_0,\ldots,\alpha_J,\beta_0,\ldots,\beta_J]$. Its {\em algebraic order} is defined to be $p$ and its {\em error constant} to be $C_{p+2}$ if, for an adequately smooth arbitrary test function $z(x)$,
\be
{\mathcal L}[h,{\bf a}]z(x)=C_{p+2}\,h^{p+2}\,z^{(p+2)}(x)+{\mathcal O}(h^{p+3}).
\label{err}
\ee
In particular,
\[
C_0=\ds\sum_{j=0}^J\alpha_j,\qquad C_1=\ds\sum_{j=0}^J j\,\alpha_j,\qquad C_q=\ds\frac{1}{q!}\ds\sum_{j=0}^Jj^q\alpha_j-\ds\frac{1}{(q-2)!}\ds\sum_{j=0}^J j^{q-2}\beta_j.
\]
The {\em principal local truncation error} ($plte$) is the leading term of (\ref{err}), i.e.,
\[
plte=C_{p+2}\,h^{p+2}\,z^{(p+2)}(x).
\]
Throughout, we shall assume that the method satisfies the following hypotheses:
\begin{enumerate}
\item $\alpha_J=1$,\hspace{1cm} $|\alpha_0|+|\beta_0|\neq 0$,\hspace{1cm}$\ds\sum_{j=0}^J|\beta_j|\neq 0$.
\item $\rho$ and $\sigma$ have no common factors.
\item $\rho(1)=\rho'(1)=0$ and $\rho''(1)=2\,\sigma(1)$; this is necessary and sufficient for the method to be consistent, that is, to have order at least one.
\item The method is zero-stable; that is, all the roots of $\rho$ lie in or on the unit circle, those on the unit circle having multiplicity not greater than one.
\end{enumerate}
The method is then convergent, see Henrici~\cite{Hen}, and the polynomial $\rho$ has root of multiplicity precisely two at $+1$.    

Algorithm~(\ref{multi}) is said to be {\em symmetric} when
\[
\alpha_j=\alpha_{J-j},\qquad \beta_j=\beta_{J-j}\qquad\mbox{for all }\, j.
\]
The algebraic order $p$ and the step number $J$ are then even numbers~\cite{LamWat}. Symmetric multistep methods are able to preserve the amplitude of the harmonic oscillator $y''=-\omega^2\,y$, see~Lambert and Watson~\cite{LamWat}. For the Schr\"odinger equation, symmetric two- and four-step methods  received particular attention in this context~\cite{Raptis1}--\cite{Simos2}. For computations of planetary orbits with symmetric multistep methods an excellent long-time behavior is reported in the literature~\cite{Quin,Ag}. A complete explanation of the behavior of classical symmetric multistep methods is given in~\cite{Hairer}. 

Usually, the coefficients of a $p$th-order linear multistep method are found from the requirement that it integrates exactly powers up to degree $p+1$, or equivalently, the operator~(\ref{func}) is vanishing for these power functions. For problems having oscillatory solutions, more efficient methods are obtained when they are exact for every linear combination of functions from the reference set
\be
\{1,x,\ldots,x^K,\exp(\pm\mu\,x),\ldots,x^P\,\exp(\pm\mu\,x)\},\qquad K+2\,P=p-1.
\label{refset}
\ee
This technique is known as {\em exponential fitting} and has a long history~\cite{Gau,Ly}. The set~(\ref{refset}) is characterized by two integer parameters, $K$ and $P$. The set in which there is no classical component is identified by $K=-1$ while the set in which there is no exponential fitting component (the classical case) is identified by $P=-1$. Parameter $P$ will be called the {\em level of tuning}. One should take in mind that exponential fitting can be applied only when a good estimate of the dominant frequency of the solution is known in advance. The coefficients of exponentially-fitted (in short: EF) methods depend on the product of the frequency $\mu$ and the stepsize $h$. An important property of EF algorithms is that they tend to the classical ones when the involved frequencies tend to zero, a fact which allows to say that exponential fitting represents a natural extension of the classical polynomial fitting. Remark that hypotheses I--IV to be convergent are only applicable to classical multistep methods. The examination of the convergence of EF multistep methods is included in Lyche's theory~\cite{Ly}. Many EF symmetric multistep methods are already constructed~\cite{Raptis1}--\cite{Simos2},~\cite{Ag}. Also, exponential fitting has been applied many times to other standard algorithms such as Runge-Kutta methods, hybrid methods, quadrature, interpolation, \ldots \,This vast material is collected by Ixaru and Vanden Berghe~\cite{Ixaru}. 

To have an idea of the accuracy of the method when solving oscillatory problems it is more appropriate to consider the {\em phase-lag}, rather than its usual $plte$. We mention the pioneering paper of Brusa and Nigro~\cite{Brus} in which the phase-lag property was introduced. This is actually another type of a truncation error, i.e. the angle between the analytical solution and the numerical solution. The purpose of this Letter is to investigate the phase properties of EF symmetric multistep methods. It turns out that for equations similar to the harmonic oscillator, the most efficient EF methods are those with the highest tuning level. In the case of the Schr\"odinger equation, this result was already obtained for particular two- and four-step EF multistep methods based on an expensive error analysis, see~\cite{Ixaru1,Ixaru2,Simos4}.
\section{Phase-lag analysis of classical symmetric multistep methods}
Linear stability analysis and phase-lag analysis of numerical methods for~(\ref{ODE2}) is based on the homogeneous test equation
\be
y'' = -\omega^2\,y,
\label{testeq}
\ee
where $\omega$ is a real constant, which may be assumed non-negative for notational convenience in latter inequalities. When we apply a symmetric multistep method to the scalar test equation~(\ref{testeq}) we obtain the difference equation
\[
\ds\sum_{j=1}^{J/2}A_j(\nu^2)\,(y_{n+j}+y_{n-j})+A_0(\nu^2)\,y_n=0,
\]
where $\nu=\omega\,h$ and $A_j(\nu^2)$ are polynomials in $\nu^2$. For a linear algorithm is $A_j(\nu^2)=\alpha_j+\nu^2\,\beta_j$, $j=0,\ldots,J/2$. Multistep methods with more than one stage give rise to higher order polynomials in $\nu^2$, see for example~\cite{Simos3} (Sec.~5.1.6).

The {\em characteristic equation} is
\be
\ds\sum_{j=1}^{J/2}A_j(\nu^2)\,(\zeta^j+\zeta^{-j})+A_0(\nu^2)=0.
\label{chareq}
\ee
Of particular interest for periodic motion is the situation where those roots are on the unit circle. A symmetric multistep method has an {\em interval of periodicity} $(0,\nu^2_0)$ if, for all $\nu\in(0,\nu^2_0)$, the roots $\zeta_i$ of the characteristic equation~(\ref{chareq}) satisfy
\be
\zeta_1=\exp(i\,\lambda(\nu)),\qquad \zeta_2=\exp(-i\,\lambda(\nu)),\qquad |\zeta_j|\leq 1,\qquad 3\leq j\leq J,
\label{percond}
\ee
where $\lambda(\nu)$ is a real function of $\nu$.
For any method corresponding to the characteristic equation~(\ref{chareq}) and for all $\nu^2\in(0,\nu_0^2)$, the {\em phase-lag} is defined as the difference
\[
t=\nu-\lambda(\nu).
\]
The {\em phase-lag order} is $q$ if
\be
t=c\,\nu^{q+1}+{\mathcal O}(\nu^{q+3}),
\label{pl}
\ee
where $c$ is the {\em phase-lag constant}.
\section{Phase-lag analysis of exponentially-fitted symmetric multistep methods}
Next, we consider EF symmetric multistep methods. In what follows, we consider the methods in their trigonometric form; i.e. instead of (\ref{refset}), the methods are exact for trigonometric functions
\[
\{1,x,\ldots,x^K,\cos(k\,x),\sin(k\,x),\ldots,x^P\,\cos(k\,x),x^P\,\sin(k\,x)\}.
\]
This is accomplished by setting $\mu=i\,k$ in the coefficients of the EF methods. The coefficients are further denoted as $\alpha_j(\theta)$ and $\beta_j(\theta)$ where $\theta=k\,h$. The corresponding classical method is determined by $\alpha_j(0)$ and $\beta_j(0)$. An application of such an EF method to the scalar test equation~(\ref{testeq}) leads to the difference equation 
\[
\ds\sum_{j=1}^{J/2}A_j(\nu^2;\theta)\,(y_{n+j}+y_{n-j})+A_0(\nu^2;\theta)\,y_n=0.
\]
The characteristic equation is
\be
\ds\sum_{j=1}^{J/2}A_j(\nu^2;\theta)\,(\zeta^j+\zeta^{-j})+A_0(\nu^2;\theta)=0,
\label{chareqef}
\ee
where $A_j(\nu^2;\theta)$ are polynomials in $\nu^2$ with coefficients which depend on the parameter $\theta$ which specifies the method of concern. For a linear algorithm is
\[
A_j(\nu^2;\theta)=\alpha_j(\theta)+\nu^2\,\beta_j(\theta),\qquad j=1,\ldots,J/2.
\]
The important modification is that the periodicity interval becomes now a two-dimensional region. Following Definition~6 of Coleman and Ixaru~\cite{ColIx}, the {\em region of stability} is a region in the $\nu-\theta$ plane, throughout which the roots of the characteristic equation~(\ref{chareqef}) satisfy the periodicity condition~(\ref{percond}). As mentioned in~\cite{Ixaru}, to study the phase properties it is more suited to replace the pair $\nu,\theta$ by the pair $\nu,r=\theta/\nu$. For any method corresponding to the characteristic equation~(\ref{chareqef}) and for all $(\nu,r\,\nu)$ belonging to the stability region, the phase-lag is defined as (\ref{pl}). The phase-lag order is $q$ if
\be
t=c(r)\,\nu^{q+1}+{\mathcal O}(\nu^{q+3}).
\label{phorder}
\ee
The following theorem was previously found by Simos and Williams~\cite{Simos3} for classical methods. Here, it is extended to the EF case. The proof is essentially the same as that given in~\cite{Simos3}. 
\begin{thm}For all $(\nu,r\,\nu)$ belonging to the stability region, a symmetric EF $J$-method with characteristic equation~(\ref{chareqef}) has phase-lag order $q$ if and only if 
\[
\ds\frac{2\,\sum_{j=1}^{J/2} A_j(\nu^2;r\,\nu)\,\cos(j\,\nu)+A_0(\nu^2;r\,\nu)}{2\,\sum_{j=1}^{J/2} j^2\,A_j(\nu^2;r\,\nu)}=-c(r)\,\nu^{q+2}+{\mathcal O}(\nu^{q+4}).
\]
\label{thm1}
\end{thm}
Our aim is to determine the expression of $c(r)$ for EF symmetric linear multistep methods.

Firstly, we reveal a relation between the $plte$ and the phase-lag. The analytical solution of (\ref{testeq}) is
\[
y(x)=c_1\,\exp(i\,\omega\,x)+c_2\,\exp(-i\,\omega\,x).
\]
Some algebraic manipulation gives
\be
{\mathcal L}[h,{\bf a}]y(x)=\Bigl(2\,\sum_{j=1}^{J/2} A_j(\nu^2;r\,\nu)\,\cos(j\,\nu)+A_0(\nu^2;r\,\nu)\Bigr)\,y(x).
\label{lte*}
\ee
Assuming that the phase-lag order is~$q$, i.e.~(\ref{phorder}), the following identity holds (see Simos and Williams~\cite{Simos3})
\be
\cos(j\,\nu)=\cos(j\,\lambda(\nu))-j^2\,c(r)\,\nu^{q+2}+{\mathcal O}(\nu^{q+4}).
\label{sim}
\ee
The roots of the characteristic equation~(\ref{chareqef}) have to satisfy the periodicity conditions~(\ref{percond}). The first two conditions of~(\ref{percond}) are equivalent to 
\be
2\,\ds\sum_{j=1}^{J/2}A_j(\nu^2;r\,\nu)\,\cos(j\,\lambda(\nu))+A_0(\nu^2;r\,\nu)=0.
\label{percondef}
\ee
We denote by $plte^*$ as the $plte$ when solving (\ref{testeq}). Using~(\ref{lte*}), (\ref{sim}) and (\ref{percondef}) we arrive to
\be
plte^*=-2\,c(r)\,\nu^{q+2}\,\ds\sum_{j=1}^{J/2} j^2\,\alpha_j(0)\,y(x).
\label{plteef1}
\ee
A convergent symmetric method has at least order two, thus $\sum_{j=1}^J j^2\,\alpha_j(0)=2\,\sum_{j=0}^J \beta_j(0)$. Note that $\sigma(1)=\sum_{j=0}^J \beta_j(0)$. By hypotheses II--III we have that $\sigma(1)\neq 0$, so $\sum_{j=0}^{J/2} j^2\,\alpha_j(0)\neq 0$. With this in mind, it follows from~(\ref{plteef1}) that the phase-lag order is equal to the order of the numerical solution of~(\ref{testeq}), a result which was already obtained by Coleman~\cite{Coleman} for classical two-step hybrid methods.

Secondly, the $plte$ of the methods considered may be written as (see~\cite{Ixaru})
\be
plte =C_{p+2}\,h^{p+2}\,(D^2+k^2)^{P+1}\,y(x),
\label{plte}
\ee
where $D^j=d^j/dx^j$ and $C_{p+2}$ is the error constant of the corresponding classical method (i.e. $\theta=0$). More specifically, when solving~(\ref{testeq}) one easily verifies that
\be
plte^*=(-1)^{p/2+1}\,C_{p+2}\,\nu^{p+2}\,(1-r^2)^{P+1}\,y(x).
\label{plteef2}
\ee
From (\ref{plte})--(\ref{plteef2}) it is clear that the order of the numerical solution of (\ref{testeq}) is equal to the order the method. Altogether, we conclude that for symmetric EF linear multistep methods, the phase-lag order is equal to the algebraic order. Very important to notice is that the algebraic order of an EF multistep method and its classical companion have both the same algebraic order~\cite{Gau,Ly,Ixaru}. Consequently, the EF procedure conserves the phase-lag order. We compare (\ref{plteef1}) and (\ref{plteef2}) to get
\be
c(r)\,=(-1)^{p/2}\,\ds\frac{C_{p+2}}{2\,\sum_{j=1}^{J/2} j^2\,\alpha_j(0)}(1-r^2)^{P+1}.
\label{cef}
\ee
Altogether, we summarize our main observations as
\begin{thm}An EF symmetric linear multistep method for (\ref{ODE2}) and the corresponding classical method have both the same phase-lag order. Furthermore, we have that
\[
c(r)=(1-r^2)^{P+1}\,c,
\]
with $r=\theta/\nu$ and $c$ is the phase-lag constant of the classical method.
\label{thmpl}
\end{thm}
When an acceptable estimate of the dominant frequency is available (i.e. $r\approx 1$) the magnitude of the phase-lag is then much smaller than that of the corresponding classical method (i.e. $r=0$). Furthermore, the more accurate the estimate of the dominant frequency is, the smaller the phase-lag is. For equations similar to (\ref{testeq}), such as the Schr\"odinger equation, it turns out that the most appropriate EF methods are those with the highest possible value of $P$. Similar results are found in~\cite{Ixaru1,Ixaru2,Simos4} for some two- and four-step EF methods for the Schr\"odinger equation via a different approach.

As an example, we select a four-step method of Simos~\cite{Simos2} (p.~351, $Case~II$) which is determined by the parameters $K=-1$ and $P=3$. Its order is six. The error constant $C_8$ of a classical symmetric four-step method reads
\[
C_8=\ds\frac{1}{20160}\,\Bigl(256+\alpha_1(0)-3584\,\beta_2(0)-56\,\beta_1(0)\Bigr).
\]
For this method is
\[
\alpha_0(0)=0,\qquad \alpha_1(0)=-1,\qquad \beta_0(0)=\ds\frac{37}{40},\qquad\qquad \beta_1(0)=\ds\frac{29}{30},\qquad \beta_2(0)=\ds\frac{17}{240},
\]
which results, with the help of~(\ref{cef}), in
\[
c(r)=\ds\frac{53}{120960}\,(1-r^2)^4.
\]
This can be easily checked with Theorem~\ref{thm1}.
\section*{Acknowledgments}
This research was supported by ``Grant 0T/04/21 of Onderzoeksfonds K.U. Leuven" and ``Scholarship BDB-B/05/06 of K.U. Leuven".

\end{document}